\documentclass[12pt]{article}
\usepackage[usenames, dvipsnames]{color}
\usepackage{authblk}
\usepackage{amsfonts}
\usepackage{mathrsfs}
\usepackage{amsmath}
\usepackage{amssymb,extarrows}
\usepackage{multicol}
\usepackage{float}
\usepackage{makeidx}
\usepackage{layout}
\usepackage{array}
\usepackage{a4wide}
\usepackage{boxedminipage}
\usepackage{hyperref}
\usepackage{latexsym}
\usepackage{color}
\usepackage{dsfont}
\usepackage{graphicx}
\usepackage{ulem}
\usepackage{amsthm}
\title{On mixed pressure-velocity regularity criteria to the Navier-Stokes equations in Lorentz spaces}
\author[1]{Hugo Beir\~{a}o da Veiga \thanks{Hugo Beir\~{a}o da Veiga (\texttt{hbeiraodaveiga@gmail.com}) partially supported  by FCT (Portugal) under the project: UIDB/MAT/04561/2020.}}
\affil[1]{\small{Department of Mathematics, Pisa University, Pisa, Italy}}
\author[2]{Jiaqi Yang\thanks{Jiaqi Yang (\texttt{yjqmath@nwpu.edu.cn}) supported by the Fundamental Research Funds for the Central Universities under grant: G2019KY05114.}}
\affil[2]{\small{School of Mathematics and Statistics, Northwestern Polytechnical University, Xi'an, 710129, China}}
\date{}

\newtheorem{theorem}{Theorem}[section]

\newtheorem{lemma}[theorem]{Lemma}

\theoremstyle{remark}
\newtheorem{remark}{Remark}[section]

\theoremstyle{definition}
\newtheorem{definition}[theorem]{Definition}

\numberwithin{equation}{section}

\newcommand{\p}{\partial}
\newcommand{\e}{\epsilon}

\newcommand{\R}{\mathbb{R}}

\newcommand{\al}{\alpha}
\newcommand{\f}{\frac}
\newcommand{\n}{\nabla}

\newcommand{\ed}{\end{document}}

\newcommand{\na}{{\nabla}}

\newcommand{\Om}{{\Omega}}
\newcommand{\Ga}{{\Gamma}}

\newcommand{\g}{{\gamma}}

\newcommand{\te}{{\theta}}

\newcommand{\m}{{\mu}}
\newcommand{\ep}{{\epsilon}}

\begin{document}
\maketitle
\begin{abstract}
In this paper we derive regular criteria in Lorentz spaces for Leray-Hopf weak solutions $v$ of the three-dimensional Navier-Stokes equations based on the formal equivalence relation $\pi\cong|v|^2$, where $\pi$ denotes the fluid pressure and $v$ the fluid velocity. It is called the mixed pressure-velocity problem (the P-V problem). It is shown that if $\f{\pi}{(e^{-|x|^2}+|v|)^{\theta}}\in L^p(0,T;L^{q,\infty})\,,$ where $0\leq\theta\leq1$ and $\f2p+\f3q=2-\theta$, then $v$ is regular on $(0,T]$. Note that, if $\Om$ is periodic, we may replace $\,e^{-|x|^2} \,$ by a positive constant. This result improves a 2018 statement obtained by one of the authors. Furthermore, as an integral part of our contribution, we give an overview on the known results on the P-V problem, and also on two main techniques used by many authors to establish sufficient conditions for regularity of the so-called Ladyzhenskaya-Prodi-Serrin (L-P-S) type.\par%
\end{abstract}

\noindent\textbf{Mathematics Subject Classification:} 35Q30, 76D03, 76D05.

\vspace{0.2cm}
\noindent \textbf{Keywords:} Navier-Stokes equations, pressure $\cong$ square velocity, regularity criteria, Lorentz spaces.

\vspace{0.2cm}

\section{Preliminaries.}
We are concerned with the regularity of weak solutions to the Navier-Stokes equations
\begin{equation}\label{NS}
\begin{cases}
\p_tv+v\cdot \n v-\Delta v+\n \pi=\,0\,,\quad&\text{in $\Omega\times(0,T)$}\,,\\
\n\cdot v=0\,,\quad&\text{in $\Omega\times(0,T)$}\,,\\
v(x,0)=v_0\,,\quad&\text{in $\Omega$}\,,
\end{cases}
\end{equation}
where the vector field $v$ is the flow velocity, the scalar function $\pi$ stands for the pressure, and the initial data $v_0$ is divergence free. In some statements an external force is assumed. Below $\,Q_T =\,\Omega \times (0,\,T]\,$, where $\,\Omega \,$ may be the whole space $\,\R^n\,;$  the $n$-dimensional torus $\,\mathbb{T}^n\,;$ or a smooth open, bounded, subset of $\,\mathbb R^n\,$. In this last case $\,\Ga\,$ denotes its boundary, and the non-slip boundary condition is always assumed:
\begin{equation}\label{nslip}
v=\,0 \quad \textrm{on} \quad \Ga \times (0,\,T]\,.
\end{equation}
Our new results concern the two first cases. The purpose of the present paper is to establish new integral criteria for regularity of solutions that relate pressure and velocity, see the left-hand side of \eqref{pcongv-teta} below. For convenience, they are called \emph{mixed pressure-velocity criteria}, abbreviate simply to \emph{P-V criteria}. It is strictly essential to start this paper by recalling  the so called Ladyzhenskaya-Prodi-Serrin (L-P-S) regularity criteria, see the pioneering references \cite{LL}, \cite{pp}, \cite{SS}. This criteria, in its \emph{strong} final form, establishes that if a weak solution $v$ of \eqref{NS} satisfies
\begin{equation}\label{LPS}
v\in L^p(0,T;L^q(\Omega))\,,\quad \f{2}{p}+\f{n}{q}=1\,,\quad q>n\,,
\end{equation}
then $\,v\,$ is a strong solution:
\begin{equation}
\label{regil}%
v \in\,L^{\infty}(0,\,T;\,H^1(\Om)\,)
\cap\,L^2(0,\,T;\,H^2(\Om)\,)\,.
\end{equation}
The result also holds for $q=n\,,$ see \cite{escau} and \cite{seregin}. Furthermore, it is well-known that strong solutions are \emph{smooth}, if data and domain are also smooth.

\vspace{0.2cm}

The long history of the condition \eqref{LPS} is completely outside of the aim of this paper. To our knowledge, the first paper where a complete proof of the above strong form was shown is Giga's 1986 reference \cite{giga}. A totally different proof was shown in the 1987 reference \cite{B87}, together with global existence results for small data, and sharp decay estimates, in the presence of general external forces. See section \ref{hbv-1987} for some details. We also recall a third distinct proof by Galdi and Maremonti in the 1988 reference \cite{GaMa}.\par%
For a "one page" proof of the L-P-S regularity criteria, in the general case \eqref{NS}, and for $\,n \geq 3\,.$ See \cite{Bv97}, by starting from equation (2.2) in this last reference.%

\vspace{0.2cm}

Coming in the wake of assumption \eqref{LPS}, many other similar sufficient conditions for regularity, but involving not just the velocity alone but also the pressure, or the gradient of the velocity, or even possible combinations, like the above P-V problem, appear. In section \ref{mot} we introduce the P-V problem and justify its relevance. Moreover, we report back on two main techniques that many authors have applied to prove the regularity criteria of the L-P-S type. Pioneering work on these two techniques, and applications to the P-V problem, are due to one of the authors, see sections \ref{restrunc} and \ref{fruitful} below. More precisely, in section \ref{restrunc} we recall and discuss the results obtained on the P-V problem by the truncation method. In section \ref{fruitful}, we recall the method introduced and developed in reference \cite{B87}, and describe the results obtained to the P-V problem by appeal to this technique.

\vspace{0.2cm}

In general, in equations like \eqref{LPS}, we put in evidence the difference between conditions with the equality sign, and conditions with this sign replaced by the inequality sign $\,<\,$. This distinction extends in an obvious way to all similar conditions considered in the sequel. For convenience, we call \emph{strong} the results in the first case, and \emph{mild-strong}, abbreviated  \emph{mild}, the results in the  second case. Weaker results are called \emph{weak}.\par%
The reader merely interested on the new results proved here may have a look to the next section and then skip directly to sections \ref{mainre}, \ref{newresults-1}. Our main result is the Theorem \ref{Thm1} below.\par%
For a rather complete introduction to the Navier-Stokes equations, from the perspective of our article, we refer to Galdi's reference \cite{Ga}.%

\section{The P-V problem and its motivation.}\label{mot}
Let's come to the main problem. The well-known equation
\begin{equation}\label{pivas}
-\Delta \pi=\sum_{i,j=1}^n\p_i\p_j(v_iv_j)\,,
\end{equation}
roughly suggests the formal equivalence
\begin{equation}\label{fr}
\pi\cong |v|^2\,.
\end{equation}
More appropriately, \eqref{pivas} merely suggests $\,{\pi \lessapprox \,|v|^2}\,,$ rather than $\, {|v|^2 \lessapprox \,\pi }\,,$ since it gives information on $\pi$ in terms of $\,v\,$, but not the reverse. This means that, formally,
$$
\frac{|\pi|}{|v|}\lessapprox\,|v|\,,
$$
but not the reverse. So results under the same integrability assumption, but on the two different quantities present in the above inequality, look stronger (more general) in the case of the left-hand side term.\par%
On the other hand, one has
$$
\frac{|\pi|}{1+\,|v|} \leq\,\frac{|\pi|}{|v|}\,.
$$
So, results obtained under conditions on the left-hand side are stronger than results under the same conditions on the right-hand side. This distinction is significant since the relation between $\pi$ and $v$ is not local. For instance, the quantity $\pi/|v|^2$ may be unbounded in some region merely due to small values of $|v|$, even if $\pi$ is bounded in the same region.\par%
The formal relation \eqref{fr} suggest the following generalization
\begin{equation}
\label{pcongv-teta}%
\frac{|\pi|}{(1+\,|v|)^\te} \cong |\,v\,|^{2-\,\te}\,.
\end{equation}
Sufficient conditions for regularity complying with \eqref{pcongv-teta} look significant since they suggest that the ties between pressure and velocity are stronger than what one could a priori expect from the global relation \eqref{pivas}. Main references on the P-V problem are \cite{B97}, \cite{B98}, \cite{B00}, \cite{zhou}, and \cite{B18}. The approach followed in the first two references, and in the three last references, are totally different. In the first couple, and for the first time, one applies De Giorgi's truncation method to the Navier-Stokes equations. This method has led to mild, instead of strong, criteria. The reason for this slight reduction of generality, actually a purely occasional fact, is quite important to the understanding of the relation between the truncation method, the functional spaces $\,L^p_*\,$ (see below), and the mild results obtained by the truncation method. This phenomena will be treated in section \ref{restrunc} below.\par%
The so-called weak-$L^p\,$ spaces, denoted in the sequel by the symbol $\,L^p_*\,,$ are just a particular case of the more recent Lorentz spaces. In fact, $\,L^p_* \equiv\,L^{p,\,\infty}\,$. See \eqref{deflorentz} below for the definition. However, in section \ref{restrunc}, we appeal to the old notation and old denomination.

\vspace{0.2cm}

The method used in references \cite{B98}, \cite{B00}, and \cite{zhou}, was developed in the 1987 reference \cite{B87}. This method has been used by many other authors, in particular by us below. A brief note on \cite{B87} will be given in section \ref{hbv-1987}. In references \cite{B98} and \cite{B00}, where $0\leq \te \leq 1\,,$ the method allows to prove strong criteria, improving the mild criteria obtained in references \cite{B97} and \cite{B98}. In reference  \cite{zhou} the case $\te>1\,$ is also treated, see below. It will be of great interest to understand why an apparent loss of regularity for the case $\te > 1$ holds.

\section{Pioneering results on the P-V criteria, and the related truncation method.}\label{restrunc}
This section mainly concerns the application of the truncation method to the P-V problem. Some words on this crucial method must be spent. The truncation method was introduced by the great mathematician Ennio De Giorgi in his outstanding 1957 paper \cite{giorgi}, where the 19th Hilbert's problem was finally solved (also solved with a different method by John Nash) after more than a half-century of attempts by many other mathematicians. The method has been further applied and developed by Guido Stampacchia in a sequence of papers, see for instance, \cite{stamp}. See also reference \cite{LUS}. Application to variational inequalities was made for the first time in 1969, by the first author of the present paper.\par%
Application to the Navier-Stokes equations, see references \cite{B97} and \cite{B98}, was in strong discontinuity with respect to the previous scalar cases. It was a considerable forward step since, in addition to the presence of a system of equations, one has also to handle the loss of the divergence-free property produced by the cut-off. Concerning the truncation method applied to the Navier-Stokes equations, we recall two very important contributions, by Vasseur \cite{vasseur} and Bjorland and Vasseur \cite{bj-va}, published respectively in 2007 and 2011. In particular, an improvement of the classical L-P-S criteria in terms of $\,L^p_w\,$ spaces is shown in \cite{bj-va}. These papers are very innovative due to the masterly use of the truncation method.

\vspace{0.2cm}

Let's turn back to the P-V problem. In \cite[Theorem 1.1]{B97} the following theorem was proved.
\begin{theorem}[see \cite{B97}, Theorem 1.1]
\label{teo-trunc-I}%
Let $\,v_0\in\,L^\infty(\Om)\cap H^1_0(\Om)\,$ be divergence free. Assume that $\,(v,\,\pi)\,$ is a weak solution to
the Navier-Stokes equations \eqref{NS} under the assumption \eqref{nslip}. Furthermore, assume that
\begin{equation}
\label{trunc-I}%
\frac{|\pi|}{1+\,|v|}\in\, L^p(0,\,T;\,L^q(\Om)\,)\,,
\end{equation}
where $\,p \in (2,\,\infty]\,,$ $\,q \in (2,\,+\infty)\,,$ and
\begin{equation}
\label{teofrac-teta}%
\frac2p +\,\frac{n}q <\,1\,.
\end{equation}
Then $v$ is bounded, and consequently is strong, and smooth if data are smooth.
\end{theorem}
Actually, assumption \eqref{trunc-I} was required merely on the subset where $\,|v|\,$ is greater than an arbitrarily large constant $\,k_0\,.$%

\vspace{0.2cm}

Furthermore, in \cite[Theorem 1.1]{B98}, the above result was extended to general values of $\theta$ with $\,0 \leq \theta \leq 1\,.$ To simplify this new attempt, it has been assumed that space and time exponents coincide, say $\,p=\,q\,.$ Following \cite{B98} we set
\begin{equation}
\label{NN}%
N=\,n+2\,.
\end{equation}
Note that $\,N\,$ is precisely the integrability exponent for which, in the particular case $\,p=\,q\,,$ the inequality \eqref{teofrac-teta} holds with the equality sign:
\begin{equation}
\label{pqN}%
\frac2N +\,\frac{n}N =\,1\,.
\end{equation}
As in \cite{B97}, the proofs given in  \cite{B98} made use of the truncation method, but with a different approach. In \cite[Theorem 1.1]{B98} the following result was proved.
\begin{theorem}[see \cite{B98}, Theorem 1.1]
\label{theo98}%
Let $\,v_0\,$ and $\,(v,\,\pi)\,$ be as in Theorem \ref{teo-trunc-I}. Furthermore, assume that for some
$\,\te \in\,[0,\,1)\,$ and for some $\,\g \in\,(2,\,N)\,,$
one has
\begin{equation}
\label{ps-98}%
\frac{|\pi|}{(1+\,|v|)^\te}\in\, L^\g_*(Q_T)\,.
\end{equation}
Then
\begin{equation}
\label{frac-98}%
v \in\, L^\m_*(Q_T\,)\,, where \quad \m=\,\m(\g)\equiv\,(1-\,\te)\,\frac{N\,\g}{N-\,\g}\,.
\end{equation}
In particular, the solution is smooth in $\,Q_T\,$ if
\begin{equation}
\label{regama-98-bis}%
\frac2\g +\,\frac n\g <\,2-\,\te\,, \quad \te \in\,[0,\,1]\,.
\end{equation}
\end{theorem}

\vspace{0.2cm}

Next we analyse the reasons which led to the mild assumption \eqref{regama-98-bis} (sign $ < $) instead of to a strong assumption (sign =). To avoid misunderstandings, note that strong assumptions lead to strong results, and mild assumptions to mild results. This is the reason for our convention, even if the strong assumption is weaker than the mild assumption. Let $\,\g_1=\,\frac{N}{2-\,\te}\,$ be the value of the parameter $\g\,$ for which \eqref{regama-98-bis} holds with the equality sign. From \eqref{frac-98} it follows that $\,\m(\g_1)=\,N\,$. So Theorem implies the following result:
\begin{equation}
\label{expl}%
\frac{|\pi|}{(1+\,|v|)^\te}\in\, L^{\g_1}_*(Q_T)\,, \quad  \frac2\g_1 +\,\frac n\g_1 =\,2-\,\te\,\quad \Longrightarrow \quad v\in L^N_*(Q_T), \quad \frac2N +\,\frac{n}N =\,1\,.
\end{equation}
Unfortunately, it is not yet know whether $  v\in L^N_*(Q_T)$ implies regularity. However, if $\g\,$ verifies \eqref{regama-98-bis}, equivalently if $\,\g >\,\g_1\,$, one has $\,\m=\,\m(\g) >\, \m(\g_1)=\,N\,.$ Since $\,L^\m_* \subset \,L^{N+\,\ep}\,$ for $\,0<\ep <\m-\,N\,,$ it follows that $\,v\in L^{N+\,\ep}(Q_T) \,,$ and smoothness follows from \eqref{LPS}. This particular case illustrates why the truncation technique has led to mild regularity statements instead of to strong statements. However, it is worth noting that the \emph{sharp} statement \eqref{expl} is not weaker than the corresponding strong statement obtained by replacing the two weak spaces in \eqref{expl} by strong Lebesgue spaces, since in the first case the right hand side (the thesis) is weaker but so is also the left hand side (the hypothesis).\par%

\vspace{0.2cm}

Let's also consider the particular case of \eqref{expl} for $\,\te=\,0\,$ (see \cite[Corollary 1.7]{B98}). For $\g_0=\,\frac N2\,,$ it follows from \eqref{frac-98} that

\begin{equation}
\label{explp}%
|\pi|\in\, L^{\g_0}_*(Q_T)\,\quad \Longrightarrow \quad v\in L^N_*(Q_T)\,.
\end{equation}
As above, to guarantee smoothness of solutions we are led to choose any given $\,\g>\,\frac N2\,.$ This leads to the following result: %
\begin{equation}
\label{explp-2}%
|\pi|\in\, L^\g_*(Q_T)\,, \quad  \frac2\g +\,\frac n\g <\,2\,\quad \Longrightarrow \quad v\in L^{N+\ep}(Q_T), \quad \frac2N +\,\frac{n}N =\,1\,,
\end{equation}
where $0<\ep<\m-\,N\,$. Hence, under the left hand side assumption, solutions are smooth (a mild regularity result).\par%
Note that regularity under P-V, or pressure alone, assumptions was simply turned into pure velocity criteria. So any improvement on velocity criteria may automatically lead to improvements on other related criteria.

\vspace{0.2cm}

A Technical Remark: In reference \cite{B98} it was assumed that $\,\g>\,2\,N/\big(2\,\te+\,(1-\,\te) \,N\big)\,.$
This assumption is superfluous \cite[Remark 1.5]{B98}. However it implies $\,\g>\,2\,,$ required in
\cite[equation (2.4)]{B98}. So, in the above formulation, this condition must be assumed.\par%

\vspace{0.2cm}

Last but not least, we refer to two interesting contributions by T. Suzuki, \cite{S1} and \cite{S2}, both in 2012, obtained by appealing to the truncation method in the \cite{B97}, \cite{B98} version. The author proves, in particular,  the following result (for details see Theorems 2.4 in \cite{S1} and 2.3 in \cite{S2}).
Assume that $p$ and $q$ satisfy \eqref{BG} below for some $\,q \in (\frac{5}{2},\,+\infty)\,.$ Then there exists $\,\ep_* >\,0\,$ such that a weak solution $\,u\,$ of the Navier-Stokes equations \eqref{NS} in $\,\R^3 \times \,(0,\,T)\,$ is smooth if it satisfies the smallness assumption
\begin{equation}
\label{trunc-IP}%
\| \pi \|_{ L^p_w(0,\,T;\,L^q_w(\Om)\,)} \leq \ep_*\,.
\end{equation}
In our context, in spite of the smallness assumption, the significance of this result is the combination of the truncation method with the condition expressed in terms of two weak $ L^p_w $ spaces.
\section{On a distinct, fruitful approach.}\label{fruitful}
The main aim of the couple of papers \cite{B00,B18} was replacing in the mixed P-V case the mild regularity assumptions by corresponding strong regularity assumptions. In the 2000 reference \cite{B00} Theorem I, the following $\theta= 1\,$ result was proved (for precise statements we always refer to the original papers).
\begin{theorem}[see \cite{B00}, Theorem I]
\label{theo2000}%
Let $\,v\,$ be a weak solution to the Navier-Stokes equations \eqref{NS} under the boundary condition \eqref{nslip}, where $\,v_0 \in L^\al(\Om) \cap \,H^1_0(\Om)\,$ is divergence free, and $\,f\in\, L^1(0,\,T;\,L^\al(\Om)\,)$ for some $\,\al>\,n\,.$  Assume that \eqref{teops-teta} and \eqref{teofrac-teta-2} below hold for $\,\te=\,1\,.$ Then
\begin{equation}
\label{reg2000}%
v \in\,C(0,\,T;\,L^\al(\Om)) \quad \textrm{and} \quad
|\,v\,|^{\,\al/2} \in L^2(0,\,T;\,H^1_0(\Om)\,)\,.
\end{equation}
In particular $v$ is smooth in $\,Q_T\,.$
\end{theorem}
The technique followed in the proof essentially appeals to the argument developed in the 1987 reference \cite{B87}. See, in particular, Lemmas 1.1 and 1.2 therein. Some information will be furnished in section \ref{hbv-1987} below.\par%
Much later, in the 2018 reference \cite{B18} Theorem 1.1, the above result was extended to the general $\,\theta $ case. Moreover the assumption $ q>n\,$ was overtook.
\begin{theorem}[see \cite{B18}, Theorem 1.1]
\label{teoaln}%
Let $\,v_0 \in L^n(\Om) \cap \,H^1_0(\Om)\,$ be divergence free and $\,f\in\, L^1(0,\,T;\,L^n(\Om)\,)$. Assume that a weak solution of the Navier-Stokes equations \eqref{NS} under the boundary condition \eqref{nslip} satisfies the assumption%
\begin{equation}
\label{teops-teta}%
\frac{|\pi|}{(1+\,|v|)^\te}\in\, L^p(0,\,T;\,L^q(\Om)\,)\,,
\end{equation}
where $\,0\leq\,\te \leq\,1\,,$ and the exponents  $\,p,\,q \in (2,\,+\infty)\,$ verify the
condition
\begin{equation}
\label{teofrac-teta-2}%
\frac2p +\,\frac{n}q =\,2-\,\te\,.
\end{equation}
If  $\,2\leq\,q <\,n\,$ we also assume that
\begin{equation}
\label{qmenosn}%
p\leq\,\frac{(n-2)\,q}{n-q}\equiv\,\frac{n-2}{(n/q)-1}\,.
\end{equation}
Under the above hypotheses one has $ \, v \in\, L^\infty(0,\,T;\,L^n(\Om)\,)\cap L^n(0,\,T;\,L^\frac{n^2}{n-2}(\Om)\,)\,,$ and $\quad \na |v|^{n/2} \in\, L^2(0,\,T;\,L^2(\Om)\,)\,.$\par%
In particular, the solution is strong. Additional smoothness of solution follows from suitable smoothness of the data.
\end{theorem}
Note that $p$ has the full range $(2,\,\infty)\,$ if $q \geq n\,.$ But for values $q<n$ the range of $p$ shrinks as $q$ decreases. For some considerations see the appendix in \cite{B18}.\par%
We advise the interested reader that notations in \cite{B87} and in \cite{B00} are different. The quantities denoted in \cite{B87} by the symbols $\,N_\al(v)\,$ and $\,M_\al(v)\,$ are the $\,\al-$powers of the quantities denoted by the same symbols $\,N_\al(v)\,$ and $\,M_\al(v)\,$ in reference
\cite{B00}. In reference \cite{B18}, see definitions (31) in this reference, the author follows the notation used in \cite{B00} for $\,\al=\,n\,.$\par%
Note that Theorem \ref{teo-trunc-I} and the last statement in Theorem \ref{theo98} are mild forms of results contained in Theorems \ref{theo2000} and \ref{teoaln} respectively. Furthermore, for $\,p=\,q=\,\g_1=\,N/(2-\,\te)\,,$ Theorem \ref{teoaln} shows that \eqref{expl} holds by replacing the two weak spaces by Lebesgue spaces.

\vspace{0.2cm}

Next, we consider the case $\,\te>\,1$  treated by Y. Zhou in the 2004 reference \cite{zhou}, by partially appealing to the method introduced in \cite{B87}. In a very systematic way, many other related criteria are proved. For the very wide set of interesting results, we refer the reader directly to the original paper. Below we will refer to the particular result concerning our main concern, namely the P-V criteria, this time for $\,\te>\,1\,.$  In this case, there is no evidence of a positive answer to the relation $\,\pi\cong |v|^2\,.$ On the contrary, both Y. Zhou's result, see below, and the constraint (51) imposed in Lemma 3.6 in reference \cite{B18}, go in the direction of a negative answer to the equivalence $\,\pi\cong |v|^2\,.$\par%
In reference \cite{zhou}, Theorem 1, item (H3), among many other results, the author states the following result (for the precise statement, see the original paper).
\begin{theorem}[see \cite{zhou}, Theorem 1, item (H3)]
\label{theozhou}%
Let $\,v_0 \in L^2(\Om) \cap \,L^q(\Om)\,,$ $\,q>\,3,$  be divergence-free, and let
$\,f=\,0\,.$  Let $\,v\,$ be a weak solution of the Navier-Stokes equations \eqref{NS} under the boundary condition \eqref{nslip}.
Furthermore, assume that $\,v\,$ satisfies \eqref{teops-teta}, where  $\,\te\in [\,1,\,5/3\,]\,,$
\begin{equation}
\label{zhou}%
\frac2p +\,\frac3q =\,\frac52 -\,\frac{3}{2}\, \te\,,
\end{equation}
and
\begin{equation}
\label{zhou-teta}%
\frac{6}{5-\,3 \te} <\,q\leq \infty\,.
\end{equation}
Then $v$ is smooth in $\,Q_T\,.$
\end{theorem}
The result extends to dimensions $\,n>\,3\,,$ see \cite{zhou}, Remark 3, item (H3)'.\par%
For the value $\te=\,1\,$ the above result coincides with the previous result obtained in reference \cite{B00}. However comparison with reference \cite{B18} looks more interesting. For $\te=\,1\,$ the two results glue perfectly. However, for $\, \te >\,1\,$ the above result looks weaker in the sense that the right hand side of \eqref{zhou} is strictly smaller than that of \eqref{fr}. Since the proofs in \cite{B00}, \cite{B18}, and \cite{zhou} have, as starting point, the ideas developed in reference \cite{B00}, we guess that all the results are the best possible attainable by the method. So the above ``loss of regularity" could be substantial, and not due to a merely technical reason. Note that larger is $\,\te\,$  "weaker" becomes the result. For $\,\te=\,\f53\,$ the assumption \eqref{zhou} becomes $\,v \in  L^\infty(Q_T)\,,$ which yields regularity by itself. It would be of \emph{great interest} having a more deep explanation of this phenomena.

\vspace{0.2cm}

Let's also consider the particular case $\te=0\,,$ the ``pressure alone" case.  A necessary classical reference is the pioneering 1969 Kaniel's paper \cite{kaniel}. In more recent times, in Berselli's reference \cite{berselli-1999}, Theorem 1.1, by following \cite{B00} (see information below), it is proved that solutions to the problem \eqref{NS}-\eqref{nslip}, satisfying the assumption
$$
\pi\in L^p(0,T;L^q(\Omega))\,, \quad \f{2}{p}+\f{n}{q}=1+ \f{n}{p}\,,\quad p<n\,,
$$
are regular. If the value $\,p=n\,$ would be reachable, the result for this particular value would be strong since the value $\,2\,$ would be attainable on the right-hand side of the above equality. Furthermore, as $p$ decreases, the result becomes weaker (for instance, for $\,p=q\,$ compare with \eqref{explp-2}).\par%
It looks useful to inform the interested readers that the item $[5]$ in the list of references in \cite{berselli-1999} was not published by the journal therein indicated. Avoiding any comment, the first author merely informs that the same paper was published, but with a different title, in another journal. It corresponds to our reference \cite{B00} below.

\vspace{0.2cm}

For the pure pressure problem in the whole space $\R^3\,$, Berselli and Galdi in reference \cite{BG}, by appealing to \cite{riogaldi} and \cite{B00}, proved regularity under the strong condition
\begin{equation}\label{BG}
\pi\in L^p(0,T;L^q(\Omega))\,, \quad \f{2}{p}+\f{n}{q}=2\,,\quad q>\f{n}{2}\,.
\end{equation}
Note that \eqref{explp-2} is a mild form of Berselli-Galdi's result. See \cite{BG} for a wide bibliography on the pressure problem.\par%
Concerning the pressure, we quote the outstanding result proved by Seregin and \v Sver\'ak, in reference \cite{sesv}. In particular the authors show that the solution is necessarily smooth if the pressure is everywhere non-negative. A result out of the main-stream, may be the more impressive global sufficient condition for regularity.\par%
To end this section we recall the 1995 reference \cite{B95} where smoothness is proved for $\,\Om=\,\R^n\,$ under assumption \eqref{BG}, this time for $\,\n v\,$ instead of $\,\pi\,$. This shows the natural equivalence between $\,\pi\,$ and $\,\n v\,.$\par%
\section{Lorentz Spaces and our Main Results.}\label{mainre}
It is worth noting that in recent years many mathematicians have been devoted to systematically extending known regularity criteria of L-P-S type from Lebesgue to Lorentz and other functional spaces. This tendency is nowadays a quite general, modern trend, in mathematics. Since Lorentz spaces are larger than Lebesgue spaces, results in Lorentz spaces are stronger than the corresponding results in Lebesgue spaces.\par%
In the 2001 reference \cite{So} H. Sohr proved that if
\begin{equation*}
v\in L^{p,s}(0,T;L^{q,\infty}(\Omega))\,,\quad \f2p+\f3q=1\,,\, 3<q<\infty\,,\, 2<s\leq p<\infty\,,
\end{equation*}
or
\begin{equation*}
\|v\|_{L^{p,\infty}(0,T;L^{q,\infty}(\Omega))}\leq\e\,,\quad \f2p+\f3q=1\,, \quad \text{for a positive constant $\e$\,,}
\end{equation*}
then the weak solution $u$ is regular on $(0,T]$. Furthermore, in the 2004 reference \cite{BM}, by appealing to the method developed in \cite{B87}, Berselli and Manfrin obtained similar results. They proved that $v$ is regular, provided that
\begin{equation*}
\|v\|_{L^{p,\infty}(0,T;L^{q}(\Omega))}\leq\e\,,\quad \f2p+\f3q=1\,,\, 3<q<\infty\,.
\end{equation*}
Recently, Suzuki in the 2012 papers \cite{S1,S2}, and Ji, Wang and Wei in the 2020 reference \cite{JWW} studied some regularity criteria in terms of the pressure $\pi$ in Lorentz spaces (we still referred to Suzuki's contributions at the end of section \ref{restrunc}, due to the appeal to the truncation method). In the 2020 reference \cite{JWW} Ji, Wang, and Wei extend Suzuki's assumption \eqref{trunc-I} to the range $\,\frac{3}{2} \leq q <\frac{5}{2}\,,$ by partially appealing to ideas in \cite{B00} and \cite{B18}.

\vspace{0.2cm}

By following the above line of research, a natural question is whether we can extend the Theorem \ref{teoaln} to Lorentz spaces. Below, we give an answer to this problem. A sufficient condition involving Lorentz spaces will be established, see equations \eqref{pisolo} and \eqref{pisolo2}. As in \cite{B18}, we may extend our new results to any space dimension $\,n \geq 3\,.$ Furthermore, extension to the boundary value problem \eqref{nslip} is the subject of a forthcoming paper.\par%
Finally, concerning some regularity criteria involving the gradient of velocity or pressure, the reader can refer to \cite{B95,BG,JWW,S1,S2}.

\vspace{0.2cm}

Let's state our new results, after recalling definition and some properties of Lorentz spaces.
\begin{definition}\label{deflorentz}
Let $1\leq p<\infty$, $1\leq q\leq \infty$, the Lorentz space $L^{p,q}$ is the set of all functions $f$ such that $\|f\|_{L^{p,q}}<\infty$, where
\begin{equation}\label{}
\|f\|_{L^{p,q}}:=
\begin{cases}
\left(p\int_0^{\infty}\tau^q|\{x\in\Omega:|f(x)|>\tau\}|^{\f{q}{p}}\f{d\tau}{\tau}\right)^{\f1q}\,,\quad& q<\infty\,,\\
\sup\limits_{\tau>0}\tau|\{x\in\Omega:|f(x)|>\tau\}|^{\f1{p}}\,,\quad& q=\infty\,.
\end{cases}
\end{equation}
\end{definition}
Actually the quantity $\,\|f\|_{L^{p,q}}\,$ is merely a semi-norm, not a norm. However it is well known that there are equivalent norms.\par%
Now, we give some useful properties which have been listed in \cite{JWW}.
\begin{description}
  \item[(i)] Interpolation character of Lorentz spaces, see for example Theorem 5.3.1 of \cite{BL},
  \begin{equation}\label{p1}
  (L^{p_0,q_0},L^{p_1,q_1})_{\delta,q}=L^{p,q}\,,\ \f{1}{p}=\f{1-\delta}{p_0}+\f{\delta}{p_1}\,,\f{1}{q}=\f{1-\delta}{q_0}+\f{\delta}{q_1}\,, 0<\delta<1\,.
  \end{equation}
  \item[(ii)] Boundedness of Riesz Transform in Lorentz spaces, see for example Lemma 2.2 of \cite{CF},
  \begin{equation}\label{p2}
  \|R_jf\|_{L^{p,q}}\leq C\|f\|_{L^{p,q}}\,,\quad 1<p<\infty\,.
  \end{equation}
  \item[(iii)] H\"{o}lder inequality in the Lorentz spaces, see for example Proposition 2.3 of \cite{LR},
  \begin{equation}\label{p3}
  \|fg\|_{L^{r,s}}\leq \|f\|_{L^{r_1,s_1}}\|g\|_{L^{r_2,s_2}}\,,
  \end{equation}
  where
  \begin{equation*}
  \f1{r}=\f1{r_1}+\f1{r_2}\,,\quad \f1s=\f1{s_1}+\f1{s_2}\,.
  \end{equation*}
  \item[(iv)] For $1\leq p<\infty$, $1\leq q_1<q_2\leq \infty$, we have, see for example Proposition 1.4.10 of \cite{G},
  \begin{equation}\label{p4}
  \|f\|_{L^{p,q_2}}\leq \left(\f{q_1}{p}\right)^{\f{1}{q_1}-\f1{q_2}}\|f\|_{L^{p,q_1}}\,.
  \end{equation}
   \item[(v)] Sobolev inequality in Lorentz spaces, see for example Theorem 8 of \cite{Tartar},
  \begin{equation}\label{p5}
  \|f\|_{L^{{\f{np}{n-p}},p}(\R^n)}\leq C\,\|\n f\|_{L^p(\R^n)}\quad \text{with $1\leq p<n$}\,.
 \end{equation}
\end{description}
Now, we state our main result.
\begin{theorem}\label{Thm1}
Set $\Omega=\text{$\R^3$ or $\mathbb{T}^3$}$. Let $(v,\pi)$ be a weak solution to \eqref{NS} with divergence-free initial data $v_0\in L^2(\Omega)\cap L^4(\Omega)$. Assume that $0\leq\theta\leq1$ and that
\begin{equation}\label{pisolo}
\f{\pi}{(e^{-|x|^2}+|v|)^{\theta}}\in L^p(0,T;L^{q,\infty}(\Omega))\,,
\end{equation}
where $p$  and  $q$ are finite, and
\begin{equation}\label{pqthe}
\f2p+\f3q =2-\theta\,.
\end{equation}
Then $v$ is regular on $(0,T] \times \Omega\,$.
\end{theorem}
\begin{remark}
When $\Omega=\mathbb{T}^3$, the assumption \eqref{pisolo} is equivalent to
\begin{equation}\label{pisolo2}
\f{\pi}{(1+|v|)^{\theta}}\in L^p(0,T;L^{q,\infty}(\Omega))\,.
\end{equation}
However, when $\Omega=\mathbb{R}^3$, the assumption \eqref{pisolo} is not replaced by \eqref{pisolo2} since we can control the term $\|(e^{-|x|^2}+|v|)^2\|^2_{L^{2}}$, but not the term $\|(1+|v|)^2\|^2_{L^{2}}$ (see equation \eqref{V-v}) according to the following proofs.\par%
We may also try to replace $1$ by a power $|v|^\m\,,$ for a suitable exponent $\m \in (0,\,1)\,,$ instead of $e^{-|x|^2}$. This would be significant, and we hope it could interest some readers.
\end{remark}
\begin{remark}
Assumption $\theta\leq1$ in our proof is necessary. In \eqref{Eq62}, the  H\"{o}lder's inequality in Lorentz spaces is used to get that
\begin{equation}\label{}
\int_{\Omega}|\tilde{\pi}|^{\beta}|\pi|^{2-\beta}V^{2+\beta\theta}dx
\leq\|\tilde{\pi}^{\beta}\|_{L^{\f{q}{\beta},\infty}}\|\pi^{2-\beta}\|_{L^{r_1,\f{2}{2-\beta}}}\|V^{2\beta}\|_{L^{r_2,\f{2}{\beta}}}.
\end{equation}
Clearly, we require $\f{2}{\beta}\geq1$. Hence $\theta=2-\f{2}{\beta}\leq1$. This constraint was already crucial in reference \cite{B18}, as explained therein. See in particular Lemma 3.6 in this last reference.
\end{remark}
\section{Proof of Theorem \ref{Thm1}}\label{newresults-1}
We first introduce the following lemma, which was proved in Lemmas 1.1 and 1.2 of \cite{B87}. See also Lemma 2.1 of \cite{B00} or Lemma 3.1 of \cite{B18}. Actually, it can be obtained by multiplying both sides of \eqref{NS} by $|v|^2v$, integrating by parts, using divergence-free condition and Cauchy-Schwarz inequality. We note that this was also the starting point of the proofs in reference \cite{JWW}.
\begin{lemma}\label{lem}
Let $(v\,,\pi)$ be a regular solution to equation \eqref{NS} in $\Omega\times[0,T]$. Then we have
\begin{equation}\label{Eq61}
\f14\f{d}{dt}\int_{\Omega}|v|^4dx+\f12\int_{\Omega}|\n v|^2|v|^2dx+\f12\int_{\Omega}|\n |v|^2|^2dx\leq \int_{\Omega}|\pi|^2|v|^2dx\,.
\end{equation}
\end{lemma}
Lemma \ref{lem} follows from the estimate (2.3) in \cite{B00} by setting $\alpha=\,4\,$ and dimension $n=3\,.$\par%

Next we set $\Omega=\R^3$ and $\beta=\f{2}{2-\theta}$, see Remark \ref{T3} for $\Omega=\mathbb{T}^3$. Note that $\beta\in[1,2]$ due to $\theta\in[0,1]$, and that $2+\theta\beta=2\beta$.\par%
Now, we control the term $\int_{\Om}|\pi|^2|v|^2dx$. For convenience, we set
\begin{equation}\label{tilde-v-pi}
V=e^{-|x|^2}+|v|\,,\quad \tilde{\pi}=\f{\pi}{(e^{-|x|^2}+|v|)^{\theta}}\,.
\end{equation}
By H\"{o}lder's inequality in Lorentz spaces \eqref{p3}, we can get
\begin{equation}\label{Eq62}
\begin{split}
\int_{\Omega}|\pi|^2|v|^2dx=&\int_{\Omega}\left(\f{\pi}{(e^{-|x|^2}+|v|)^{\theta}}\right)^{\beta}|\pi|^{2-\beta}(e^{-|x|^2}+|v|)^{\beta\theta}|v|^{2}dx\\
\leq&\int_{\Omega}|\tilde{\pi}|^{\beta}|\pi|^{2-\beta}V^{2+\beta\theta}dx\\
\leq&\|\tilde{\pi}^{\beta}\|_{L^{\f{q}{\beta},\infty}}\|\pi^{2-\beta}\|_{L^{r_1,\f{2}{2-\beta}}}\|V^{2\beta}\|_{L^{r_2,\f{2}{\beta}}}\\
=&\|\tilde{\pi}\|^{\beta}_{L^{q,\infty}}\|\pi\|^{2-\beta}_{L^{(2-\beta)r_1,2}}
\|V^{2}\|^{\beta}_{L^{\beta r_2,2}}\,,
\end{split}
\end{equation}
where
\begin{equation}\label{relation1}
\f{\beta}{q}+\f1{r_1}+\f{1}{r_2}=1\,.
\end{equation}
Here, we remark that when $\theta=1$, i.e. $\beta=2$, the corresponding estimate is
\begin{equation}\label{}
\int_{\Omega}|\pi|^2|v|^2dx\leq\|\tilde{\pi}\|^{2}_{L^{q,\infty}}\|V^{2}\|^{2}_{L^{\f{2q}{q-2},2}}\,.
\end{equation}
Since $-\Delta\pi=\sum^3_{i,j=1}\p_i\p_j(v_iv_j)$, by \eqref{p2}, we have
\begin{equation}\label{65}
\begin{split}
\int_{\Omega}|\pi|^2|v|^2dx\leq& C\|\tilde{\pi}\|^{\beta}_{L^{q,\infty}}\|\pi\|^{2-\beta}_{L^{(2-\beta)r_1,2}}
\|V^{2}\|^{\beta}_{L^{\beta r_2,2}}\\
\leq& C\|\tilde{\pi}\|^{\beta}_{L^{q,\infty}}\||v|^2\|^{2-\beta}_{L^{(2-\beta)r_1,2}}
\|V^{2}\|^{\beta}_{L^{\beta r_2,2}}\\
\leq&C\|\tilde{\pi}\|^{\beta}_{L^{q,\infty}}\|V^2\|^{2-\beta}_{L^{(2-\beta)r_1,2}}
\|V^{2}\|^{\beta}_{L^{\beta r_2,2}}\,.
\end{split}
\end{equation}
By the interpolation character of Lorentz spaces \eqref{p1} and by Sobolev inequality in Lorentz spaces \eqref{p5}, it follows that
\begin{equation}\label{abcd}
\|V^2\|_{L^{(2-\beta)r_1,2}}\leq C\|V^2\|^{1-\delta_1}_{L^{2,2}}
\|V^2\|^{\delta_1}_{L^{6,2}}\leq C \|V^2\|^{1-\delta_1}_{L^{2}}\|\n V^2\|^{\delta_1}_{L^{2}}
\end{equation}
and
\begin{equation}\label{efg}
\|V^2\|_{L^{\beta r_2,2}}\leq
C\|V^2\|^{1-\delta_2}_{L^{2,2}}\|V^2\|^{\delta_2}_{L^{6,2}}
\leq C\|V^2\|^{1-\delta_2}_{L^{2}}\|\n V^2\|^{\delta_2}_{L^{2}}\,,
\end{equation}
where $0<\delta_1\,,\delta_2<1$, and
\begin{equation}\label{relation2}
\f{1}{(2-\beta)r_1}=\f{1-\delta_1}{2}+\f{\delta_1}{6}\,,\quad\f{1}{\beta r_2}=\f{1-\delta_2}{2}+\f{\delta_2}{6}\,.
\end{equation}
We remark that there exist $r_1$ and $r_2$ satisfying \eqref{relation1} and \eqref{relation2}. Actually, we can take
\[
\f{1}{r_1}=\f{2-\beta}{2}\left(1-\f{\beta}{q}\right)\,,\quad \f{1}{r_2}=\f{\beta}{2}\left(1-\f{\beta}{q}\right)\,,
\]
and therefore $\delta_1=\delta_2=\f{3\beta}{2q}\in(0,1)$ due to $q\in(\f{3}{2-\theta},\infty)$.
Hence, from \eqref{65} it follows that
\begin{equation}\label{}
\begin{split}
\int_{\Omega}|\pi|^2V^2dx \leq&C\|\tilde{\pi}\|^{\beta}_{L^{q,\infty}}\|V^2\|^{(1-\delta_1)(2-\beta)}_{L^{2}}
\|\ \n V^2\|^{\delta_1(2-\beta)}_{L^2}
\|V^2\|^{(1-\delta_2)\beta}_{L^{2}}\|\n V^2\|^{\delta_2\beta}_{L^2}\\
\leq&C\|\tilde{\pi}\|^{\beta}_{L^{q,\infty}}\|V^2\|^{(1-\delta_1)(2-\beta)+(1-\delta_2)\beta}_{L^{2}}
\|\n V^2\|^{\delta_1(2-\beta)+\delta_2\beta}_{L^{2}}\\
\leq&C\|\tilde{\pi}\|^{\f{2\beta}{2-\delta_1(2-\beta)-\delta_2\beta}}_{L^{q,\infty}}
\|V^2\|^{\f{2[(1-\delta_1)(2-\beta)+(1-\delta_2)\beta]}{2-\delta_1(2-\beta)-\delta_2\beta}}_{L^{2}}
+\e\|\n V^2\|^{2}_{L^{2}}\,.
\end{split}
\end{equation}
Noting that
\begin{equation}\label{}
(1-\delta_1)(2-\beta)+(1-\delta_2)\beta
=2-\delta_1(2-\beta)-\delta_2\beta\,,
\end{equation}
we have
\[
\f{2[(1-\delta_1)(2-\beta)+(1-\delta_2)\beta]}{2-\delta_1(2-\beta)-\delta_2\beta}=2\,.
\]
Thus we have
\begin{equation}\label{seisonze}
\begin{split}
\int_{\Omega}|\pi|^2|v|^2dx
\leq&C_\e \,\|\tilde{\pi}\|^{\f{2\beta}{2-\delta_1(2-\beta)-\delta_2\beta}}_{L^{q,\infty}}
\|V^2\|^{2}_{L^{2}}
+\e\|\n V^2\|^{2}_{L^{2}}\,.
\end{split}
\end{equation}
Note that
\begin{equation}\label{V-v}
\|V^2\|^{2}_{L^{2}}=\|e^{-2|x|^2}+2e^{-|x|^2}|v|+|v|^2\|^2_{L^{2}}\leq  C(1+\|v\|^2_{L^{2}}+\||v|^2\|^2_{L^{2}})
\end{equation}
and
\begin{equation}\label{nV-v}
\begin{split}
\|\n V^2\|^{2}_{L^{2}}=&\|\n (e^{-2|x|^2}+2e^{-|x|^2}|v|+|v|^2)\|^2_{L^{2}}\\
\leq& C(1+\|v\|^2_{L^{2}}+\|\n v\|^2_{L^{2}}+\|\n|v|^2\|_{L^{2}})\,.
\end{split}
\end{equation}
Hence, we have
\begin{equation}\label{seisonze2}
\begin{split}
\int_{\Omega}|\pi|^2|v|^2dx
\leq&C_\e \,\|\tilde{\pi}\|^{\f{2\beta}{2-\delta_1(2-\beta)-\delta_2\beta}}_{L^{q,\infty}}(1+\|v\|^2_{L^{2}}
+\||v|^2\|^2_{L^{2}})\\
&+C\e(1+\|v\|^2_{L^{2}}+\|\n v\|^2_{L^{2}})+C\e\|\n |v|^2\|^{2}_{L^{2}}\,.
\end{split}
\end{equation}
By this estimate and Lemma \ref{lem}, and setting $\e$ sufficiently small, using Gronwall's lemma and Ladyzhenskaya-Prodi-Serrin regularity criteria \eqref{LPS}, we can get that $v$ is smooth in $\Omega\times[0,T]\,,$ provided that
\[
\tilde{\pi}\in L^{\f{2\beta}{2-\delta_1(2-\beta)-\delta_2\beta}}(0,T;L^{q,\infty})\,.
\]
Finally, if we have
\begin{equation}\label{}
2\f{2-\delta_1(2-\beta)-\delta_2\beta}{2\beta}+\f3{q}=2-\theta\,,
\end{equation}
then we can get Theorem \ref{Thm1}. Actually, from \eqref{relation1} and \eqref{relation2}, we have
\begin{equation}\label{}
\begin{split}
&1=\f{\beta}{q}+\f1{r_1}+\f{1}{r_2}\\
&=\f{\beta}{q}+(2-\beta)\left(\f{1-\delta_1}{2}+\f{\delta_1}{6}\right)+\beta\left(\f{1-\delta_2}{2}+\f{\delta_2}{6}\right)\\
&=\f{\beta}{q}+\f12(2-\beta)+\f12\beta-\f13\delta_1(2-\beta)-\f13\delta_2\beta\\
&=\f{\beta}{q}+1-\f13\delta_1(2-\beta)-\f13\delta_2\beta\,,
\end{split}
\end{equation}
which gives
\begin{equation}\label{}
\begin{split}
\delta_1(2-\beta)+\delta_2\beta
=\f{3\beta}{q}\,.
\end{split}
\end{equation}
Hence, we have
\begin{equation}\label{}
\begin{split}
2\,\f{2-\delta_1(2-\beta)-\delta_2\beta}{2\beta}+\f3{q}
=\f{2}{\beta}-\f{\delta_1(2-\beta)+\delta_2\beta}{\beta}+\f3{q}=\f2{\beta}=2-\theta\,.
\end{split}
\end{equation}
\begin{remark}\label{T3}
When $\Omega=\mathbb{T}^3$, the Sobolev inequality in Lorentz spaces should be
\begin{equation}
\|f-\int_{\mathbb{T}^3}fdx\|_{L^{{\f{np}{n-p}},p}(\mathbb{T}^n)}\leq C\,\|\n f\|_{L^p(\mathbb{T}^n)}\quad \text{with $1\leq p<n$}\,.
\end{equation}
Hence,
\begin{equation}
\|f\|_{L^{{\f{np}{n-p}},p}(\mathbb{T}^n)}\leq C\,(\|f\|_{L^p(\mathbb{T}^n)}+\|\n f\|_{L^p(\mathbb{T}^n)})\quad \text{with $1\leq p<n$}\,.
\end{equation}
Thus, as the main differences of the proofs, the above \eqref{abcd} and \eqref{efg} should be replaced by
\begin{equation}\label{abcd2}
\|V^2\|_{L^{(2-\beta)r_1,2}}\leq C\|V^2\|^{1-\delta_1}_{L^{2,2}}
\|V^2\|^{\delta_1}_{L^{6,2}}\leq C \|V^2\|^{1-\delta_1}_{L^{2}}(\|V^2\|+\|\n V^2\|_{L^{2}})^{\delta_1}
\end{equation}
and
\begin{equation}\label{efg2}
\|V^2\|_{L^{\beta r_2,2}}\leq
C\|V^2\|^{1-\delta_2}_{L^{2,2}}\|V^2\|^{\delta_2}_{L^{6,2}}
\leq C\|V^2\|^{1-\delta_2}_{L^{2}}(\|V^2\|+\|\n V^2\|_{L^{2}})^{\delta_2}\,,
\end{equation}
respectively, and therefore \eqref{seisonze} becomes
\begin{equation}\label{seisonze-2}
\begin{split}
\int_{\Omega}|\pi|^2|v|^2dx
\leq&C_\e \,\|\tilde{\pi}\|^{\f{2\beta}{2-\delta_1(2-\beta)-\delta_2\beta}}_{L^{q,\infty}}
\|V^2\|^{2}_{L^{2}}+\e\|V^2\|^{2}_{L^{2}}
+\e\|\n V^2\|^{2}_{L^{2}}\,.
\end{split}
\end{equation}
Remaining proofs are the same.
\end{remark}
\section{Notes on reference \cite{B87}.}\label{hbv-1987}
By taking into account that some ideas developed in reference \cite{B87} has been a main departure point in the proofs of many of the results quoted in the previous sections, it looks suitable to give here some comments (due to the first author) on the above publication, first published as the IMA preprint \cite{B85}.\par%
To our knowledge, a complete proof of the L-P-S \emph{strong} condition for regularity \eqref{LPS} was shown for the first time on Y.Giga's 1986 reference \cite{giga}. A totally different proof was also shown in the 1987 reference \cite{B87}, even if this fact was not explicitly written as a formal theorem (see below). Reference \cite{B87} was received for publication on October 1985, hence without intersection with the 1986 reference \cite{giga} (received for publication on July 1984). A third distinct proof was given by G.P.Galdi and P.Maremomti in the 1988 reference \cite{GaMa}.%

\vspace{0.2cm}

The proof in reference \cite{B87} was given up to obvious details, already well known at that time. We briefly explain this point below where, for convenience, we replace the $(q,\alpha)-\,$notation used in \cite{B87} by our present notation $\,(p,\,q)\,.$ The following is one of the results proved in \cite{B87}.

\begin{theorem}[see \cite{B87}, Theorem 0.1]
\label{teo87}%
Consider the evolution Navier-Stokes equations in the whole space $\,\R^n\,$, with a divergence free initial data $\,v_0 \in L^q(\R^n)\,$ and an external force $\,f \in \,L^1(0,T\, L^q(\R^n)\,.$ Assume that
\begin{equation}
\label{AA}%
v\in\, L^p(0,\,T;\,L^q\,)\,, \quad \,q>\,n\,,
\end{equation}
where
\begin{equation}
\label{BB}%
\frac2p +\,\frac{n}q =\,1\,.
\end{equation}
Under the above hypothesis, if $v$ is "sufficient regular", one has
\begin{equation}\label{CC}%
\,\|v(t)\|_q \leq\, exp \,\Big( c \mu^{-(n+q)/(q-n)} \|v\|^p_{L^p(0,t;L^q)}\,\Big)\,\big(\,\|v_0\|_q +\, \|f\|_{L^1(0,t;L^q)}\,\big)
\end{equation}
for every $\,t \in [0,\,T]\,.$ In particular $\,v \in\, L^\infty(0,\,T;\,L^q(\R^n))\,.$
\end{theorem}
Note that \eqref{AA} under assumption \eqref{BB} coincides with assumption \eqref{LPS}. Moreover, the estimate \eqref{CC} implies, in particular, $\,v \in\, L^\infty(0,\,T;\,L^q(\R^n))\,.$ This immediately yields smoothness of solutions, due (for instance) to a previous 1983 result by H. Sohr, see \cite{sohr}, or much more simpler, by appealing to weaker (or mild) forms of assumption \eqref{LPS}, which were well know and discussed at that time. This situation was claimed in \cite{B87}, Remark (i), on page 152.\par%
As remarked in \cite{B87} page 153, in the proof of Theorem 0.1 the author proves a priori estimates in the usual mathematical sense. To justify formal calculations, some additional regularity on the solution was assumed. Obviously, this regularity was not used to estimate any kind of quantities.\par%
The following bold type remark was stated immediately after Theorem 0.1 in \cite{B87}: \emph{``The a priori estimate \eqref{CC} can be utilized to show that if a solution $v$ of (0.1) belongs to the class $\,L^p(0,\,T; L^q)\,,$ then $\,v\in L^\infty(0,\,T;\,L^q)\,,$ and \eqref{CC} holds"}.\par%
 Furthermore, this sentence was immediately followed by this second remark: "We leave the technical details to the interested reader. Note that the existence of a solution in the class $\, L^p(0,\,T;\,L^q\,)\,$ is an open problem".%

\vspace{0.2cm}

The fact that the existence of a solution in the above class \eqref{AA}-\eqref{BB} was an open problem led the author, at that time, to avoid an explicit statement (a theorem) merely based on a conjecture. In fact, full $\,C^\infty(Q_T)\,$ regularity was explicitly stated as Theorems each time the additional L-P-S assumption was not required in the proof. This was the case for the results under smallness assumptions on initial data and external forces like, for instance, global $\,C^\infty(Q_\infty)\,$ regularity for sufficient small data. See \cite{B87}, Theorems 0.2 and 0.3, and Remark (i), page 152. See also Theorems 2.1 and 2.2. In these cases, non-additional conditions of regularity were assumed, and this allowed explicit formal theorems.\par%
In \cite{B87}, many local and global sharp estimates were also proved, in particular, lower and upper bounds on time, and decay at infinity.%

\begin{thebibliography}{99}
\bibitem{B85}
H. Beir\~{a}o da Veiga, \textit{Existence and asymptotic behaviour for strong solutions of the Navier-Stokes equations in the whole space,} IMA Preprint Series, preprint n. 190, October 1985. Institute for Mathematics and its Applications, University of Minnesota, Minneapolis.
%
\bibitem{B87}
H. Beir\~{a}o da Veiga, \textit{Existence and asymptotic behaviour for strong solutions of the Navier-Stokes equations in the whole space,} Indiana Univ. Math. J., \textbf{36} (1987), 149--166.
\bibitem{B95}
H. Beir\~{a}o da Veiga, \textit{A new regularity class for the Navier-Stokes equations in $R^n$,} Chin. Ann. Math. Ser. B \textbf{16} (1995), 407--412.
\bibitem{Bv97}
H. Beir\~{a}o da Veiga, \textit{Remarks on the smoothness of the $\,L^{\infty}(0,\,T; L^3\,)\,$ solutions of the $\,3-D\,$ Navier-Stokes equations,} Portugaliae Math., \textbf{54} (1997), 381--391.
\bibitem{B97}
H. Beir\~{a}o da Veiga, \textit{Concerning the regularity of the solutions to the Navier-Stokes equations via the truncation method; Part I,} Diff. Int. Eq., \textbf{10} (1997), 1149--1156.
\bibitem{B98}
H. Beir\~{a}o da Veiga, \textit{Concerning the regularity of the solutions to the Navier-Stokes equations via the truncation method. Part II,} \'{E}quations aux D\'{e}riv\'{e}es Partielles et Applications; Articles d\'{e}di\'{e}s \`{a} J.L. Lions \`{a} l'occasion de son 70. anniversaire; Gauthier-Villars, Paris (1998), 127--138.
\bibitem{B00}
H. Beir\~{a}o da Veiga, \textit{A sufficient condition on the pressure for the regularity of weak solutions to the Navier-Stokes equations,} J.Math. Fluid Mech. \textbf{2} (2000), 99--106.
\bibitem{B18}
H. Beir\~{a}o da Veiga, \textit{On the Truth, and Limits, of a Full Equivalence $p \cong v^2$ in the Regularity Theory of the Navier-Stokes Equations: A Point of View}, J.Math. Fluid Mech. \textbf{20} (2018), 889--898.
\bibitem{BL}
J. Bergh,, J. L\"{o}fstr\"{o}m, \textit{Interpolation Spaces}, Springer, Berlin, 1976.
%
\bibitem{berselli-1999}
L.C.~Berselli, \textit{sufficient conditions for the regularity of the solutions of the
Navier--Stokes equations}, Math. Meth. Appl. Sci., \textbf{22} (1999), 1079--1085.
%
\bibitem{BG}
L.C. Berselli, G.P. Galdi, \textit{Regularity criteria involving the pressure for the weak solutions to the Navier-Stokes equations,} Proc. Am. Math. Soc. \textbf{130} (2002), 3585--3595.
\bibitem{BM}
L.C. Berselli and R. Manfrin, \textit{On a theorem of Sohr for the Navier-Stokes equations,} J. Evol. Eq., \textbf{4} (2004), 193--211.
\bibitem{bj-va}
C. Bjorland, A.F.Vasseur, \textit{Weak in space, Log in time improvement of the Lady\u zenskaja-Prodi-Serrin criteria,} J. Math. Fluid Mech., \textbf{13} (2011), 259--269.
%
\bibitem{CF}
J.A. Carrillo, L.C.F. Ferreira, \textit{Self-similar solutions and large time asymptotics for the dissipative quasi-geostrophic
equation}, Monatsh. Math. \textbf{151} (2007), 111--142.
%
\bibitem{escau}
L.~Escauriaza, G.~Seregin, and V.~\v Sver\'ak, \emph{$\,L_{3,\,\infty}\,$
solutions to the Navier--Stokes equations and backward uniqueness},
Russian Math. Surveys \textbf{58} (2003), 211--250.
\bibitem{giorgi}
E. De Giorgi,  \textit{Sulla differenziabilit\`a e l'analicit\`a delle estremali degli integrali multipli regolari}, Mem. Accad. Sci. Torino, cl. Sci. Fis. Mat. Nat. (3), 3 (1957), 25--43.
%
\bibitem{Ga}
G.P. Galdi, \textit{An introduction to the Navier-Stokes initial-boundary value problem,} Fundamental Directions in Mathematical Fluid-Mechanics, Birkhauser, Basel, 2000, p.p. 1-70. MR 2002c:35207.
%
\bibitem{GaMa}
G.P. Galdi, P. Maremonti,  \textit{Sulla regolarit\`a delle soluzioni deboli al sistema di Navier-Stokes in domini arbitrari,} Ann. Univ. Ferrara., \textbf{34} (1988), 59--73.
\bibitem{giga}
Y. Giga,  \textit{Solutions for semilinear parabolic equations in $\,L^p\,$ and regularity of weak solutions of the Navier-Stokes system }, J. Diff. Eq., \textbf{61} (1986), 186--212.
\bibitem{G}
L. Grafakos, \textit{Classical Fourier Analysis,} 2nd edn, Springer, Berlin, 2008.
\bibitem{JWW}
X. Ji, Y. Wang, W. Wei, \textit{New regularity criteria based on pressure or gradient of velocity in Lorentz spaces for the 3D Navier-Stokes equations,} J.Math. Fluid Mech., \textbf{22} (2020), 1--8.
\bibitem{kaniel}
S.~Kaniel, \textit{A sufficient condition for smoothness of solutions of Navier--Stokes equations}, Israel J.Math., \textbf{6} (1969), 354--358.
%
\bibitem{LL}
O.A.~Lady\v zhenskaya, \textit{Uniqueness and smoothness of generalized solutions of Navier--Stokes equations}, Zap. Nau\v cn. Sem. Leningrad Otdel. Mat. Inst. Steklov (LOMI), \textbf{5} (1967) 169--185.
%
\bibitem{LUS}
O.A.~Lady\v zhenskaya, N.N. Ural'ceva, and V.A. Solonnikov, \textit{Linear and Quasilinear Equations of Parabolic Type}, Amer. Math. Soc., Providence, R.I., 1968 (translated from Russian).
%
\bibitem{LR}
P.G. Lemari\'{e}-Rieusset, \textit{Recent developments in the Navier-Stokes problem,} Chapman \& Hall/CRC, London.
%
\bibitem{pp}
G.~Prodi, \textit{Un teorema di unicit\`a per le equazioni di Navier--Stokes}, Ann. Mat. Pura Appl., \textbf{48} (1959), 173--182.
%
\bibitem{riogaldi}
S.~Rionero and G.P. Galdi, \textit{The weight function approach to uniquiness of viscous flows in unbounded domains}, Arch. Rat. Mech. Anal., \textbf{69} (1979), 37--52.
%
\bibitem{seregin}
G.~Seregin, \emph{On smoothness of $\,L_{3,\,\infty}-$ solutions to the Navier--Stokes equations up to the boundary}, Math. Ann.,
\textbf{332} (2005), 219--238.
%
\bibitem{sesv}
G.~Seregin and V.~\v Sver\'ak, \textit{Navier-Stokes equations with lower bounds on the pressure}, Arch. Rat. Mech. Anal.\textbf{163} (2002),
65--86.
%
\bibitem{SS}
J. Serrin, \textit{The initial value problem for the Navier-Stokes equations}, in R.E. Langer editor, Nonlinear Problems, Univ. Wisconsin Press, Madison, Wisconsin (1963), 69--98.
%
\bibitem{sohr}
H.~Sohr, \textit{Zur Regularit\"atstheorie der instation\"aren Gleichungen von Navier-Stokes},\, Math. Z.\,\textbf{184}\, (1983),
359--375.
\bibitem{So}
H. Sohr, \textit{A regularity class for the Navier-Stokes equations in Lorentz spaces,} J. Evol. Equ. \textbf{1} (2001), 441--467.
%
\bibitem{stamp}
G. Stampacchia,  \textit{Le probl\`eme de Dirichlet pour les \' equations elliptiques du second ordre a coefficients discontinus}, Ann. Inst. Fourier Grenoble,  \textbf{15} (1965), 189--258.

\bibitem{S1}
T. Suzuki, \textit{Regularity criteria of weak solutions in terms of the pressure in Lorentz spaces to the Navier-Stokes equations,}
J. Math. Fluid Mech. \textbf{14} (2012), 653--660.
%
\bibitem{S2}
T. Suzuki, \textit{A remark on the regularity of weak solutions to the Navier-Stokes equations in terms of the pressure in
Lorentz spaces,} Nonlinear Anal. Theory Methods Appl. \textbf{75} (2012), 3849--3853.
\bibitem{Tartar}
L. Tartar, \textit{Imbedding theorems of Sobolev spaces into Lorentz spaces}, Boll. dell'Unione Mat. Ital. \textbf{1} (1998), 479--500.
\bibitem{vasseur}
A.F. Vasseur, \textit{A new proof of partial regularity of solutions to Navier-Stokes equations,} NoDEA., \textbf{14} (2007), 753--785.
\bibitem{zhou}
Y.~Zhou, \textit{Regularity criteria in terms of pressure for the $3$-$D$  Navier-Stokes equations}, Math. Ann., \textbf{328} (2004), 173--192.
%
\end{thebibliography}
\end{document}